\documentclass[leqno,12pt]{article} 

\usepackage{amsmath}
\usepackage{amssymb}
\usepackage{pst-all}

\newtheorem{theorem}{Theorem}[section]

\newtheorem{definition}[theorem]{Definition}
\newtheorem{example}[theorem]{Example}

\newtheorem{proposition}[theorem]{Proposition}
\newtheorem{corollary}[theorem]{Corollary}

\newtheorem{lemma}[theorem]{Lemma}

\newtheorem{conjecture}{Conjecture}[section]

\newcommand{\R}{\mathbb{R}}
\newcommand{\N}{\mathbb{N}}

\newcommand{\qed}{\ensuremath{\quad \blacksquare} \vspace{\baselineskip}}

\newcommand{\Diag}{\mathrm{Diag\,}}
\newcommand{\diag}{\mathrm{diag\,}}

\begin{document}
\allowdisplaybreaks
\title{\textbf{Generalized Hadamard Product and the Derivatives of Spectral Functions}}
\author{\textbf{
Hristo S.\ Sendov\thanks{Department of Mathematics and Statistics,
University of Guelph, Guelph, Ontario, Canada N1G 2W1. Email:
\texttt{hssendov\char64 uoguelph.ca}. Research supported by
NSERC.}
 }}

\maketitle

\begin{abstract}
In this work we propose a generalization of the Hadamard product
between two matrices to a tensor-valued, multi-linear product
between $k$ matrices for any $k \ge 1$. A multi-linear dual
operator to the generalized Hadamard product is presented. It is a
natural generalization of the $\Diag x$ operator, that maps a
vector $x \in \R^n$ into the diagonal matrix with $x$ on its main
diagonal. Defining an action of the $n \times n$ orthogonal
matrices on the space of $k$-dimensional tensors, we investigate
its interactions with the generalized Hadamard product and its
dual. The research is motivated, as illustrated throughout the
paper, by the apparent suitability of this language to describe
the higher-order derivatives of spectral functions and the tools
needed to compute them. For more on the later we refer the reader
to \cite{Sendov2003b} and \cite{Sendov2003c}, where we use the
language and properties developed here to study the higher-order
derivatives of spectral functions.
\end{abstract}



\section{Introduction}

Spectral functions, are functions on a symmetric matrix argument
invariant under a closed subgroup of the orthogonal group on the
space of all $n \times n$ symmetric matrices, $S^n$. More
precisely, $F : S^n \rightarrow \R$ is spectral if
$$
F(U^TXU) = F(X),
$$
for all $X \in S^n$ and $U \in O(n)$ --- the orthogonal group on
$\R^n$. It is not difficult to see that such functions can be
represented as the composition
$$
F = f \circ \lambda,
$$
where $f : \R^n \rightarrow \R$ is a symmetric function
($f(Px)=f(x)$ for any permutation matrix $P$ and vector $x$), and
$\lambda : S^n \rightarrow \R^n$ is the eigenvalue map:
$\lambda(X) = (\lambda_1(X),...,\lambda_n(X))$ --- all eigenvalues
of $X$. We will assume throughout that,
$$
\lambda_1(X) \ge \cdots \ge \lambda_n(X).
$$

The study of spectral functions generalizes the study of the
individual eigenvalues of a symmetric matrix since if we let
\begin{align*}
\phi_k(x) &: \R^n \rightarrow \R, \\
\phi_k(x) &:= \mbox{ the $k^{\mbox{\scriptsize th}}$ largest
element of } \{x_1,...,x_n \},
\end{align*}
then $\phi_k(x)$ is symmetric and
$$
\lambda_k(X) = (\phi_k \circ \lambda) (X).
$$

Various smoothness properties of eigenvalues have been studied for
some time now and find a lot of applications in areas ranging from
matrix perturbation theory \cite{StewardSun:1990}, and eigenvalue
optimization \cite{LewisOverton:1996}, \cite{Lewis:2003}, to
quantum mechanics \cite{kemble:1958}. The Taylor expansion (when
it exists) of the eigenvalues of symmetric matrices depending on
one scalar parameter are described in the monograph by Kato
\cite{kato:1976full}. This naturally raises the questions about
the differentiability properties of the spectral functions. Many
such questions have already been investigated in the literature
(see below) and surprisingly the answers to most of them follow
the same pattern: $f \circ \lambda$ has a property at the matrix
$X$ if, and only if, $f$ has the same property at the vector
$\lambda(X)$. It is only natural, then to try to describe the
differentials of $f \circ \lambda$ in terms of the differentials
of the simpler function $f$.

Here is a list of properties for which $f \circ \lambda$ has that
property at (or around) the matrix $X$ if and only if $f$ has the
same property at (or around) the vector $\lambda(X)$.
\begin{enumerate}
\item $F$ is lower semicontinuous at $A$ if, and only if, $f$ is
at $\lambda(A)$, \cite{Lewis:1994a}. \item $F$ is lower
semicontinuous and convex if, and only if, $f$ is,
\cite{Davis:1957}, \cite{Lewis:1994a}. \item The symmetric
function corresponding to the Fenchel conjugate of $F$ is the
Fenchel conjugate of $f$, \cite{seeger:1997}, \cite{Lewis:1994a}.
(A similar statement holds for the recession function of $F$,
\cite{seeger:1997}.) \item  $F$ is pointed, has good asymptotic
behaviour or is a barrier function on the set $\lambda^{-1}(C)$
if, and only if, $f$ is such on $C$, \cite{seeger:1997}. \item $F$
is Lipschitz around $A$ if, and only if, $f$ is such around
$\lambda(A)$, \cite{Lewis:1994b} \item \label{diff:list} $F$ is
(continuously) differentiable at $A$ if, and only if, $f$ is at
$\lambda(A)$, \cite{Lewis:1994b}. \item $F$ is strictly
differentiable at $A$ if, and only if, $f$ is at $\lambda(A)$,
\cite{Lewis:1994b}, \cite{Lewis:1996}. \item $\nabla(f \circ
\lambda)$ is semismooth at $X$ if, and only if, $\nabla f$ is at
$\lambda(X)$, \cite{QiYang:2003}. \item If $f$ is l.s.c. and
convex, then $F$ is twice epi-differentiable at $A$ relatively to
$\Omega$ if, and only if, $f$ is twice epi-differentiable at
$\lambda(A)$ relative to $\lambda(\Omega)$,
\cite{TorkiM:Jan301997c}, where $\Omega$ is an arbitrary
epi-gradient. \item $F$ has a quadratic expansion at $X$ if, and
only if, $f$ has a quadratic expansion at $\lambda(X)$,
\cite{LewisSendov:2000}. \item $F$ is twice (continuously)
differentiable at $X$ if, and only if, $f$ is twice (continuously)
differentiable at $\lambda(X)$, \cite{LewisSendov:2000a}. \item
\label{infty::list} $F \in \mathcal{C}^{\infty}$ at $A$
$\Leftrightarrow f \in \mathcal{C}^{\infty}$ at $\lambda(A)$,
\cite{Dadok:1982}. \item \label{analitic:list} $F$ is analytic at
$A$ if, and only if, $f$ is at $\lambda(A)$,
\cite{TsingFanVerriest:1994}. \item $F$ is a polynomial of the
entries of $A$ if, and only if, $f$ is a polynomial.  This is a
consequence of the Chevalley Restriction Theorem,
\cite[p.~143]{Warner:1972}.
\end{enumerate}

We want to stress that there are exceptions to the pattern. For
example if $f$ is directionally differentiable at $\lambda(X)$
this doesn't imply that $f \circ \lambda$ is such at $X$, see
\cite{Lewis:1994b}.

In \cite{Lewis:1994b} and \cite{LewisSendov:2000a} the authors
gave explicit formulae for the gradient and the Hessian of the
spectral function $F$ in terms of the derivatives of the symmetric
function $f$. In order to reproduce them here we need a bit more
notation. For any vector $x$ in $\R^n$, Diag$\,x$ will denote the
diagonal matrix with vector $x$ on the main diagonal, and diag$\,:
M^n \rightarrow \R^n$ will denote its dual operator defined by
diag$\,(X)=(x_{11},...,x_{nn})$. Recall that the Hadamard product
of two matrices $A=[A^{ij}]$ and $B=[B^{ij}]$ of the same size is
the matrix of their element-wise product $A \circ B =
[A^{ij}B^{ij}]$. Thus we have
\begin{align}
\label{deriv:1} \nabla (f \circ \lambda)(X) &= V \big(\Diag \nabla
f(\lambda(X)) \big) V^T, \mbox{ and } \\
\label{deriv:2} \nabla^2 (f \circ \lambda)(X)[H_1,H_2] &= \nabla^2
f(\lambda(X))[\mbox{\rm diag}\, \tilde{H_1}, \mbox{\rm diag}\,
\tilde{H_2}] + \\
\nonumber & \hspace{4cm} + \langle \mathcal{A}(\lambda(X)),
\tilde{H_1} \circ \tilde{H_2} \rangle,
\end{align}
where $V$ is any orthogonal matrix such that $X=V \big(\mbox{\rm
Diag}\, \lambda(X) \big) V^T$ is the {\it ordered} spectral
decomposition of $X$; $\tilde{H_i} = V^TH_iV$ for $i=1,2$, and
$x\in \R^n \rightarrow \mathcal{A}(x)$ is a matrix valued map that
is continuous if $\nabla^2 f(x)$ is.

In \cite{LewisSendov:2000a} a conjecture was made that $F$ is
$k$-times (continuously) differentiable at $A$ if, and only if,
$f$ is such at $\lambda(A)$. It is conceivable that high-powered
analytical methods may give a direct proof of this conjecture, but
never the less an interesting question is what the
$k^{\mbox{\scriptsize th}}$ differential of $F$ looks like and how
to compute it practically. Explicit formula for the
$k^{\mbox{\scriptsize th}}$ differential of $F$ will generalize
the formula for the $k^{\mbox{\scriptsize th}}$ term in the Taylor
expansion (when it exists) of the individual eigenvalues given in
\cite{kato:1976full}.

Before attacking the questions in the previous paragraph we need
to answer several more basic questions. What are the common
features in Formulae (\ref{deriv:1}) and (\ref{deriv:2}), that we
expect to generalize when we further differentiate?  We propose a
language that shows a good promise to simplify the description of
the higher order derivatives of spectral functions. It is based on
the idea of generalizing the Hadamard product of two matrices to a
$k$-tensor valued product between $k$ matrices. The current paper
is the first of three. It defines what we mean by a generalized
Hadamard product and investigates some of its multi-linear
algebraic properties. In \cite{Sendov2003b} we will formulate
calculus-type rules for the interaction between the generalized
Hadamard product and the eigenvalues of symmetric matrices.
Finally, in \cite{Sendov2003c} we will describe how to compute the
derivatives of spectral functions in two important cases. In
particular, we will show that Conjecture~\ref{conj2} holds for the
derivatives of any spectral function at a symmetric matrix with
distinct eigenvalues, as well as for the derivatives of separable
spectral functions at an arbitrary symmetric matrix. (Separable
spectral functions are those arising from symmetric functions
$f(x) = g(x_1) + \cdots + g(x_n)$ for some function $g$ on a
scalar argument.)


\section{Generalizations of the Hadamard product}
\label{GenHadProd}

By $\{H_{pq} : 1 \le p, q \le n \}$ we will denote the standard
basis of the space of all $n \times n$ matrices. That is, the
matrices $H_{pq}$ are such that $(H_{pq})^{ij}$ is $1$ if
$(i,j)=(p,q)$, and $0$ otherwise.

Let us look closely at the Hadamard product, $H_1 \circ H_2$,
between two matrices $H_1$ and $H_2$ from $M^n$. It is a matrix
valued function on two matrix arguments, linear in each argument
separately. Therefore it is uniquely determined by its values on
the pairs of basic matrices $(H_{p_1q_1}, H_{p_2q_2})$.

On such basic pairs the Hadamard product is defined as:
$$
(H_{p_1q_1} \circ H_{p_2q_2})^{ij} = \left\{
\begin{array}{ll}
1, & \mbox{ if } i=p_1=p_2 \mbox{ and } j=q_1=q_2, \\
0, & \mbox{ otherwise. }
\end{array}
\right.
$$
Naturally, we may define the {\it cross Hadamard product} by the
rule
$$
(H_{p_1q_1} \circ_{\hspace{-0.05cm}_{\mbox{\tiny $(12)$}}}
H_{p_2q_2})^{ij} := \left\{
\begin{array}{ll}
1, & \mbox{ if } i=p_1=q_2 \mbox{ and } j=p_2=q_1, \\
0, & \mbox{ otherwise, }
\end{array}
\right.
$$
and then extend this to a bilinear function on all $M^n \times
M^n$. The Hadamard product and the cross Hadamard product are
essentially the same thing:
$$
H_{p_1q_1} \circ_{\hspace{-0.05cm}_{\mbox{\tiny $(12)$}}}
H_{p_2q_2} = H_{p_1q_1} \circ H_{p_2q_2}^T = H_{p_1q_1} \circ
H_{q_2p_2}.
$$

These observations can be naturally generalized in the following
way. Denote the set $\{1,2,...,k\}$ of the first $k$ natural
numbers by $\N_k$. A {\it $k$-tensor on $\R^n$} is a real-valued
map on $\R^n \times \cdots \times \R^n$ ($k$-times) linear in each
argument separately. When a basis in $\R^n$ is fixed, a $k$-tensor
can be viewed as an $n \times \cdots \times n$ ($k$-times)
``block'' of numbers. We will index the elements of a tensor just
like we index the entries of a matrix. The space of all
$k$-tensors on $\R^n$ will be denoted by $T^{k,n}$.

\begin{definition} \rm
\label{tripleHad} For a fixed permutation $\sigma$ on $\N_k$, we
define {\it $\sigma$-Hadamard} product between $k$ matrices to be
a $k$-tensor on $\R^n$ as follows. Given any $k$ basic matrices
$H_{p_1q_1}$, $H_{p_2q_2}$,...,$H_{p_kq_k}$ we define:
\begin{align*}
 (H_{p_1q_1} \circ_{\sigma} H_{p_2q_2}
\circ_{\sigma} \cdots \circ_{\sigma} H_{p_kq_k})^{i_1i_2...i_k} &=
\left\{
\begin{array}{ll}
1, \hspace{-0.2cm} & \mbox{ if } i_s=p_s=q_{\sigma(s)}, \forall s=1,...,k, \\
0, \hspace{-0.2cm} & \mbox{ otherwise. }
\end{array}
\right.
\end{align*}
 Now, extend this product to a $k$-tensor valued map on $k$ matrix arguments,
linear in each of them separately.
\end{definition}

Another way to write the above definition is using the Kronecker
delta. Recall that $\delta_{ij}$ is equal to $1$ if $i=j$, and $0$
otherwise. Thus,
\begin{align}
\label{boz} \nonumber (H_{p_1q_1} \circ_{\sigma} H_{p_2q_2}
\circ_{\sigma} \cdots \circ_{\sigma} H_{p_kq_k})^{i_1i_2...i_k} &=
\delta_{i_1p_1}\delta_{i_1q_{\sigma(1)}}\cdots
\delta_{i_{k}p_{k}}\delta_{i_{k}q_{\sigma(k)}} \\[-0.3cm]
& \\[-0.3cm]
& \nonumber = \delta_{i_1p_1}\delta_{p_1q_{\sigma(1)}}\cdots
\delta_{i_{k}p_{k}}\delta_{p_{k}q_{\sigma(k)}}.
\end{align}

The next lemma gives the formula for the general entry of the
$\sigma$-Hadamard product between arbitrary matrices.

\begin{lemma}
\label{Hadam:arb} The $\sigma$-Hadamard product of arbitrary
matrices is given by
\begin{align*}
(H_{1}\circ_{\sigma} H_{2} \circ_{\sigma} \cdots \circ_{\sigma}
H_k)^{i_1i_2...i_k} &= H_1^{i_1i_{\sigma^{-1}(1)}}\cdots
H_k^{i_ki_{\sigma^{-1}(k)}} \\
&=H_{\sigma(1)}^{i_{\sigma(1)}i_1}\cdots
H_{\sigma(k)}^{i_{\sigma(k)}i_k}.
\end{align*}
\end{lemma}

\begin{proof}
Let $\sigma$ be a permutation on $\N_k$ and let $H_1$,...,$H_k$ be
arbitrary matrices. Using the definition that the product is
linear in each argument separately, we compute
\begin{align*}
(H_{1} &\circ_{\sigma} H_{2} \circ_{\sigma} \cdots \circ_{\sigma}
H_k)^{i_1i_2...i_k} \\ &= \sum_{p_1,q_1 = 1}^{n,n} \cdots
\sum_{p_k,q_k = 1}^{n,n} H_1^{p_1q_1}\cdots
H_k^{p_kq_k}(H_{p_1q_1} \circ_{\sigma} H_{p_2q_2} \circ_{\sigma}
\cdots
\circ_{\sigma} H_{p_{k}q_{k}})^{i_1i_2...i_k} \\
&=\sum_{p_1, q_1 = 1}^{n,n} \cdots \sum_{p_k, q_k = 1}^{n,n}
H_1^{p_1q_1}\cdots H_k^{p_kq_k}
\delta_{i_1p_1}\delta_{i_1q_{\sigma(1)}}\cdots
\delta_{i_{k}p_{k}}\delta_{i_{k}q_{\sigma(k)}} \\
&=\sum_{p_1, q_1 = 1}^{n,n} \cdots \sum_{p_k, q_k = 1}^{n,n}
H_1^{p_1q_1}\cdots H_k^{p_kq_k}
\delta_{i_1p_1}\delta_{i_{\sigma^{-1}(1)}q_{1}}\cdots
\delta_{i_{k}p_{k}}\delta_{i_{\sigma^{-1}(k)}q_{k}} \\
&=H_1^{i_1i_{\sigma^{-1}(1)}}\cdots H_k^{i_ki_{\sigma^{-1}(k)}} \\
&=H_{\sigma(1)}^{i_{\sigma(1)}i_1}\cdots
H_{\sigma(k)}^{i_{\sigma(k)}i_k}. \\[-1.6cm]
\end{align*}
\hfill \qed
\end{proof}

\begin{corollary}
When the first $k-1$ of the matrices involved in the product are
basic we get
\begin{align*}
&(H_{p_1q_1} \circ_{\sigma} \cdots \circ_{\sigma}
H_{p_{k-1}q_{k-1}} \circ_{\sigma} H)^{i_1i_2...i_k} \\
&=\delta_{i_1p_1} \delta_{i_1 q_{\sigma(1)}} \cdots
\delta_{i_{l-1}p_{l-1}} \delta_{i_{l-1}q_{\sigma(l-1)}}
H^{i_{\sigma(l)}i_{l}} \delta_{i_{l+1}p_{l+1}}
\delta_{i_{l+1}q_{\sigma(l+1)}} \cdots
\delta_{i_{k}p_{k}}\delta_{i_{k}q_{\sigma(k)}},
\end{align*}
where $l=\sigma^{-1}(k)$.
\end{corollary}

\begin{proof}
Let $l=\sigma^{-1}(k)$, using the result of the previous lemma we
calculate.
\begin{align*}
&(H_{p_1q_1} \circ_{\sigma} \cdots \circ_{\sigma}
H_{p_{k-1}q_{k-1}} \circ_{\sigma} H)^{i_1i_2...i_k}
=H_{p_1q_1}^{i_1i_{\sigma^{-1}(1)}}\cdots
H_{p_{k-1}q_{k-1}}^{i_{k-1}i_{\sigma^{-1}(k-1)}}
H^{i_ki_{\sigma^{-1}(k)}} \\
&= \delta_{i_1p_1}\delta_{i_{\sigma^{-1}(1)}q_1}\cdots
\delta_{i_{k-1}p_{k-1}}\delta_{i_{\sigma^{-1}(k-1)}q_{k-1}}
H^{i_ki_{\sigma^{-1}(k)}} \\
&=\delta_{i_1p_1} \delta_{i_1 q_{\sigma(1)}} \cdots
\delta_{i_{l-1}p_{l-1}} \delta_{i_{l-1}q_{\sigma(l-1)}}
H^{i_{\sigma(l)}i_{l}} \delta_{i_{l+1}p_{l+1}}
\delta_{i_{l+1}q_{\sigma(l+1)}} \cdots
\delta_{i_{k}p_{k}}\delta_{i_{k}q_{\sigma(k)}}. \\[-1.4cm]
\end{align*}
\hfill \qed
\end{proof}

The above corollary can be easily modified when the matrix $H$ is
in arbitrary position in the product.

\begin{example} \rm
We already saw that, when $k=2$ and $\sigma = (12)$ the
$\sigma$-Hadamard product is essentially the ordinary Hadamard
product:
$$
H_1 \circ_{\hspace{-0.05cm}_{\mbox{\tiny $(12)$}}} H_2 = H_1 \circ
H_2^T.
$$
If we restrict our attention to the space of symmetric matrices,
then the two products coincide.  In the case when $\sigma =
(1)(2)$ we get
$$
H_1 \circ_{\hspace{-0.05cm}_{\mbox{\tiny $(1)(2)$}}} H_2 = (\diag
H_1)(\diag H_2)^T.
$$
\end{example}

\begin{example} \rm
The $\sigma$-Hadamard product has meaning even when $k=1$. In that
case, there is just one permutation on the set $\N_1$ and the
$\sigma$-Hadamard product corresponding to it has one matrix
argument and returns, by definition, a vector ($1$-tensor). Since
$\sigma = (1)$, extending the notation, the $\sigma$-Hadamard
product is given by the rule:
\begin{align*}
(\circ_{\sigma} H_{p_1q_1})^{i_1} &=  \left\{
\begin{array}{ll}
1, & \mbox{ if } i_1=p_1=q_1 \\
0, & \mbox{ otherwise }
\end{array}
\right.\\
&= (\diag H_{q_1p_1})^{i_1}.
\end{align*}
Extending by linearity we get
$$
\circ_{\sigma} H = \diag H.
$$
\end{example}

For any two $k$-tensors, $T_1$, and $T_2$ we define a scalar
product between them in the natural way:
$$
\langle T_1, T_2 \rangle = \sum_{i_1,...,i_k =1}^{n,...,n}
T_1^{i_1...i_k} T_2^{i_1...i_k}.
$$

\begin{lemma}
\label{dot-prod} Let $T$ be a $k$-tensor on $\R^n$, and $H$ be a
matrix in $M^n$. Let $H_{p_1q_1}$,...,$H_{p_{k-1}q_{k-1}}$ be
basic matrices in $M^n$, and let $\sigma$ be a permutation on
$\N_k$. Then the following identities hold.
\begin{enumerate}
\item If $\sigma^{-1}(k)=k$, then
\begin{align*}
\langle T, H_{p_1q_1} &\circ_{\sigma} \cdots \circ_{\sigma}
H_{p_{k-1}q_{k-1}} \circ_{\sigma} H \rangle =
\Big(\prod_{t=1}^{k-1} \delta_{p_tq_{\sigma(t)}} \Big)
\sum_{t=1}^n T^{p_1 \dots p_{k-1}t}H^{tt}.
\end{align*}
\item If $\sigma^{-1}(k) = l$, where $l \not= k$, then
\begin{align*}
\langle T, H_{p_1q_1} &\circ_{\sigma} \cdots \circ_{\sigma}
H_{p_{k-1}q_{k-1}} \circ_{\sigma} H \rangle =
\Big(\prod_{\substack{t=1 \\ t \not= l}}^{k-1}
\delta_{p_tq_{\sigma(t)}} \Big) T^{p_1 \dots
p_{k-1}q_{\sigma(k)}}H^{q_{\sigma(k)}p_{\sigma^{-1}(k)}}.
\end{align*}
\end{enumerate}
\end{lemma}

\begin{proof}
Using the definitions and observation (\ref{boz}), we calculate.
\begin{align*}
&\langle T, H_{p_1q_1} \circ_{\sigma} H_{p_2q_2} \circ_{\sigma}
\cdots \circ_{\sigma} H_{p_{k-1}q_{k-1}} \circ_{\sigma} H \rangle\\
&=\sum_{p_k, q_k = 1}^{n,n} H^{p_kq_k} \langle T, H_{p_1q_1}
\circ_{\sigma} H_{p_2q_2} \circ_{\sigma} \cdots \circ_{\sigma}
H_{p_{k-1}q_{k-1}} \circ_{\sigma} H_{p_kq_k} \rangle \\
&= \sum_{p_k, q_k = 1}^{n,n} H^{p_kq_k}
\sum_{i_1,...,i_k=1}^{n,...,n} T^{i_1...i_k} (H_{p_1q_1}
\circ_{\sigma} H_{p_2q_2} \circ_{\sigma} \cdots \circ_{\sigma}
H_{p_{k-1}q_{k-1}} \circ_{\sigma} H_{p_kq_k})^{i_1...i_k} \\
&=\sum_{p_k, q_k = 1}^{n,n} H^{p_kq_k}
\sum_{i_1,...,i_k=1}^{n,...,n} T^{i_1...i_k}
\delta_{i_1p_1}\delta_{p_1q_{\sigma(1)}}\cdots
\delta_{i_{k}p_{k}}\delta_{p_{k}q_{\sigma(k)}} \\
&=\sum_{p_k, q_k = 1}^{n,n} H^{p_kq_k} T^{p_1...p_k}
\delta_{p_1q_{\sigma(1)}} \cdots \delta_{p_{k}q_{\sigma(k)}}.
\end{align*}
The result follows easily by considering the two cases separately.
\hfill \qed
\end{proof}

\section{A partial order on $P^k$ and one property of the $\sigma$-Hadamard product}
\label{sect:part-ord}

Given two permutations $\sigma$, $\mu$ on $\N_k$, we say that
$\sigma$ {\it refines} $\mu$ if for every $s \in \N_k$ there is an
$r \in \N_k$ such that
$$
\{\sigma^l(s) \, : \, l=1,2,... \} \subseteq \{\mu^l(r) \, : \,
l=1,2,... \},
$$
where $\sigma^l(s)= \sigma ( \sigma ( \cdots ( \sigma(s))\cdots )$
- $l$ times.

In other words $\sigma$ {\it refines} $\mu$ if every cycle of
$\sigma$ is contained in a cycle of $\mu$. Clearly the cycles of
$\sigma$ will partition the cycles of $\mu$. If $\sigma$ refines
$\mu$ we will denote it by
$$
\mu \preceq \sigma.
$$
The set of all permutations on $\N_k$ as well as the set of all $n
\times n$ permutation matrices will be denoted by $P^k$. Clearly
the refinement is a pre-order on $P^k$ (it is reflexive,
transitive, but not antisymmetric). With respect to this
pre-order, the identity permutation is the biggest element (that
is, bigger that any one else) and every permutation with only one
cycle is a smallest element (that is, it is smaller than any other
element).

There is a natural map between the set $P^k$ and the {\it diagonal
subspaces} of $\R^k$, given as follows:
\begin{align*}
\mathcal{D}(\sigma) 
&=\{ x \in \R^k \, : \, x_{s} = x_{\sigma(s)} \,\, \forall s \in
\N_k \}.
\end{align*}
This map is onto but is not one-to-one since, for example, when
$k=3$ $\mathcal{D}((123)) = \mathcal{D}((132))=\{x \in \R^3 \, :
\, x_1=x_2=x_3\}$. Clearly the image of the identity permutation
is $\R^k$. The following relationship helps to visualize the
partial order on $P^k$
$$
\mu \preceq \sigma \,\,\, \Leftrightarrow  \,\,\, \mathcal{D}(\mu)
\subseteq \mathcal{D}(\sigma).
$$

Finally, given a tensor $T \in T^{k,n}$ we may want to preserve
the entries lying on a diagonal ``subspace'' of $T$ and substitute
the rest of the entries of $T$ with zeros. In other words, given a
permutation $\mu \in P^k$, we introduce the notation $P_{\mu}(T)$
for the tensor in $T^{k,n}$ defined by
$$
(P_{\mu}(T))^{i_1...i_k} = \left\{ \begin{array}{ll}
T^{i_1...i_k}, & \mbox{ if } i_{s} = i_{\mu(s)}, \,  \forall s \in \N_k \\
0, & \mbox{ otherwise. }
\end{array}
\right.
$$

After all these preparations, we can formulate the main result in
this section. It describes when we can transfer diagonal
``subspaces'' of $T$ between different $\sigma$-Hadamard products.

\begin{theorem}
\label{trans:thm} Let $\sigma_1$, $\sigma_2$, and $\mu$ be three
permutations on $\N_k$. Then the identity
\begin{equation*}
\langle P_{\mu}(T), H_1 \circ_{\sigma_1} \cdots \circ_{\sigma_1}
H_k \rangle = \langle P_{\mu}(T), H_1 \circ_{\sigma_2} \cdots
\circ_{\sigma_2} H_k \rangle
\end{equation*}
holds for any matrices $H_1$,...,$H_k$, and any tensor $T$ in
$T^{k,n}$ if, and only if, $\mu \preceq \sigma_2^{-1} \circ
\sigma_1$.
\end{theorem}

\begin{proof}
Since both sides are linear in each of the matrices $H_1$,..,$H_k$
separately, it is enough to prove the theorem when these matrices
are basic. In other words, we are going to show that
\begin{equation*}
\langle P_{\mu}(T), H_{p_1q_1} \circ_{\sigma_1} \cdots
\circ_{\sigma_1} H_{p_kq_k} \rangle = \langle P_{\mu}(T),
H_{p_1q_1} \circ_{\sigma_2} \cdots \circ_{\sigma_2} H_{p_kq_k}
\rangle,
\end{equation*}
for any indexes $p_1$,...,$p_k$, $q_1$,...,$q_k$, and for any $T
\in T^{k,n}$ if, and only if, $\mu \preceq \sigma_2^{-1} \circ
\sigma_1$.  Direct calculation shows:
\begin{align*}
\langle P_{\mu}(T), H_{p_1q_1} &\circ_{\sigma_1} \cdots
\circ_{\sigma_1} H_{p_kq_k} \rangle \\
&= \sum_{i_1,...,i_k=1}^{n,...,n} (P_{\mu}(T))^{i_1...i_k}
(H_{p_1q_1} \circ_{\sigma_1} \cdots \circ_{\sigma_1}
H_{p_kq_k})^{i_1...i_k} \\
&=\sum_{i_1,...,i_k=1}^{n,...,n} (P_{\mu}(T))^{i_1...i_k}
H_{p_1q_1}^{i_1i_{\sigma_1^{-1}(1)}} \cdots
H_{p_kq_k}^{i_ki_{\sigma_1^{-1}(k)}} \\
&=\sum_{i_1,...,i_k=1}^{n,...,n}
(P_{\mu}(T))^{i_1...i_k}\delta_{i_1p_1}
\delta_{i_1q_{\sigma_1(1)}}
\cdots \delta_{i_kp_k} \delta_{i_kq_{\sigma_1(k)}} \\
&=(P_{\mu}(T))^{p_1...p_k} \delta_{p_1q_{\sigma_1(1)}} \cdots
\delta_{p_kq_{\sigma_1(k)}}.
\end{align*}
The last expression is equal to $T^{p_1...p_k}$ when
$p_s=p_{\mu(s)}=q_{\sigma_1(s)}$ for all $s \in \N_k$, and is
equal to $0$ otherwise.

Analogously we have
\begin{align*}
\langle P_{\mu}(T), H_{p_1q_1} \circ_{\sigma_2} \cdots
\circ_{\sigma_2} H_{p_kq_k} \rangle &= (P_{\mu}(T))^{p_1...p_k}
\delta_{p_1q_{\sigma_2(1)}} \cdots \delta_{p_kq_{\sigma_2(k)}},
\end{align*}
which is equal to $T^{p_1...p_k}$ when
$p_s=p_{\mu(s)}=q_{\sigma_2(s)}$ for all $s \in \N_k$, and is
equal to $0$ otherwise.

Suppose that $\mu \preceq \sigma_2^{-1} \circ \sigma_1$.  We
consider three cases.

If there is an $s_0$ such that $p_{s_0} \not= p_{\mu(s_0)}$, then
both expressions are zero and we trivially have equality.

If $p_s = p_{\mu(s)}$ for all $s \in \N_k$ but for some $s_0$ we
have that $p_{s_0} \not= q_{\sigma_1(s_0)}$, then it is not
possible to have $p_s=q_{\sigma_2(s)}$ for all $s \in \N_k$.
Indeed, suppose on the contrary that $p_s=q_{\sigma_2(s)}$ for all
$s \in \N_k$. Letting $r=\sigma_2(s)$ we get $p_{\sigma_2^{-1}(r)}
= q_r$ for every $r \in \N_k$. Therefore
$p_{\sigma_2^{-1}(\sigma_1(s))} = q_{\sigma_1(s)}$ for every $s
\in \N_k$. In particular $p_{\sigma_2^{-1}(\sigma_1(s_0))} =
q_{\sigma_1(s_0)} \not= p_{s_0}$. But $\mu \preceq \sigma_2^{-1}
\circ \sigma_1$ implies that $\sigma_2^{-1}(\sigma_1(s_0))$ and
$s_0$ belong to the same cycle of $\mu$, that is $\mu^{l}(s_0) =
\sigma_2^{-1}(\sigma_1(s_0))$ for some $l \in \N$. By the
assumption in this case we have that $p_{s_0}=p_{\mu^l(s_0)}$ for
every $l$. This is a contradiction. Thus for some $s_1 \in \N_k$
we have $p_{s_1} \not= q_{\sigma_2(s_1)}$ and again we will have
that both expressions are equal to zero.

Suppose finally that $p_s=p_{\mu(s)}=q_{\sigma_1(s)}$ for all $s
\in \N_k$. Then the first expression is equal to $T^{p_1...p_k}$.
If we show that $p_s=q_{\sigma_2(s)}$ for every $s \in \N_k$, then
we will be done. Suppose this is not true, that is, for some
$s_0$, $p_{s_0} \not= q_{\sigma_2(s_0)}$. Then for $r_0 =
\sigma_2(s_0)$ we will have $p_{\sigma_2^{-1}(r_0)} \not=
q_{r_0}$, and for $s_1 = \sigma_1^{-1}(r_0)$ we have
$p_{\sigma_2^{-1}(\sigma_1(s_1))} \not= q_{\sigma_1(s_1)}$. Again
$\mu \preceq \sigma_2^{-1} \circ \sigma_1$ implies that
$\sigma_2^{-1}(\sigma_1(s_1))$ and $s_1$ belong to the same cycle
of $\mu$ and we reach a contradiction as in the previous case.

To prove the opposite direction of the theorem, suppose that
\begin{equation}
\label{jenkris} (P_{\mu}(T))^{p_1...p_k}
\delta_{p_1q_{\sigma_1(1)}} \cdots \delta_{p_kq_{\sigma_1(k)}} =
(P_{\mu}(T))^{p_1...p_k} \delta_{p_1q_{\sigma_2(1)}} \cdots
\delta_{p_kq_{\sigma_2(k)}},
\end{equation}
for every choice of the indexes $p_1$,...,$p_k$ and
$q_1$,...,$q_k$ and every $T$. Take $T$ to be such that
$T^{i_1...i_k} \not= 0$ for every choice of the indexes
$i_1,...,i_k$ satisfying $i_s=i_{\mu(s)}$ for every $s \in \N_k$.
Suppose that $\mu \npreceq \sigma_2^{-1} \circ \sigma_1$. This
means that there is an number $s_0 \in \N_k$ such that
$\sigma_2^{-1}( \sigma_1(s_0))$ and $s_0$ are not in the same
cycle of $\mu$. Choose the indexes $p_1$,...,$p_k$ and
$q_1$,...,$q_k$ so that $p_s = p_{\mu(s)}$ and
$p_s=q_{\sigma_1(s)}$, for every $s \in \N_k$. Moreover, choose
the indexes $p_1$,...,$p_k$ so that if $s,r \in \N_k$ are not in
the same cycle of $\mu$, then $p_s \not= p_r$. This in particular
means that $p_{\sigma_2^{-1}( \sigma_1(s_0))} \not= p_{s_0}$.

With the choices so made, the left-hand side of
Equation~(\ref{jenkris}) will be equal to $T^{i_1...i_k} \not= 0$.
We will reach a contradiction if we show that for some $r_0$,
$p_{r_0} \not= q_{\sigma_2(r_0)}$, since then the right-hand side
of Equation~(\ref{jenkris}) will be zero. Suppose on the contrary
that $p_{r} = q_{\sigma_2(r)}$ for every $r \in \N_k$. Then,
$$
p_{\sigma_2^{-1}(\sigma_1(s))} = q_{\sigma_1(s)} = p_s, \,\,\,
\mbox{ for every } s \in \N_k.
$$
Substitute above $s=s_0$ to reach a contradiction. Thus, $p_{r_0}
\not= q_{\sigma_2(r_0)}$ for some $r_0 \in N_k$ and we are done.
\hfill \qed
\end{proof}

Notice that if $\mu \preceq \nu$, then for arbitrary permutation
$\sigma$ in $P^k$ we have
$$
\mu \preceq \nu^{-1} = (\sigma \circ \nu)^{-1} \circ \sigma.
$$
This observation leads to the next corollary.

\begin{corollary}
Suppose $\mu$ and $\nu$ are permutations in $P^n$ such that $\mu
\preceq \nu$. Then for arbitrary permutation $\sigma \in P^k$, any
matrices $H_1$,...,$H_k$, and a tensor $T$ in $T^{k,n}$ we have
the identity:
\begin{equation*}
\langle P_{\mu}(T), H_1 \circ_{\sigma} \cdots \circ_{\sigma} H_k
\rangle = \langle P_{\mu}(T), H_1 \circ_{\sigma \circ \nu} \cdots
\circ_{\sigma \circ \nu} H_k \rangle.
\end{equation*}
In particular, the result holds when $\nu = \mu$ or $\nu =
\mu^{-1}$.
\end{corollary}

It will be useful to see what are the conclusions of the above
theorem when $k \le 3$. We summarize them in the next corollary.

\begin{corollary}
\label{ident:cor} For any $T \in T^{2,n}$ and any two matrices
$H_1$ and $H_2$ we have
$$
\langle P_{\hspace{-0.05cm}_{\mbox{\tiny $(12)$}}}(T), H_1
\circ_{\hspace{-0.05cm}_{\mbox{\tiny $(1)(2)$}}} H_2 \rangle =
\langle P_{\hspace{-0.05cm}_{\mbox{\tiny $(12)$}}}(T), H_1
\circ_{\hspace{-0.05cm}_{\mbox{\tiny $(12)$}}} H_2 \rangle.
$$
For any $T \in T^{3,n}$ and any three matrices $H_1$, $H_2$, and
$H_3$ we have
\begin{align*}
\langle P_{\hspace{-0.05cm}_{\mbox{\tiny $(13)$}}}(T), H_1
\circ_{\hspace{-0.05cm}_{\mbox{\tiny $(132)$}}} H_2
\circ_{\hspace{-0.05cm}_{\mbox{\tiny $(132)$}}} H_3 \rangle &=
\langle P_{\hspace{-0.05cm}_{\mbox{\tiny $(13)$}}}(T), H_1
\circ_{\hspace{-0.05cm}_{\mbox{\tiny $(12)(3)$}}} H_2
\circ_{\hspace{-0.05cm}_{\mbox{\tiny $(12)(3)$}}} H_3 \rangle, \\
\langle P_{\hspace{-0.05cm}_{\mbox{\tiny $(23)$}}}(T), H_1
\circ_{\hspace{-0.05cm}_{\mbox{\tiny $(123)$}}} H_2
\circ_{\hspace{-0.05cm}_{\mbox{\tiny $(123)$}}} H_3 \rangle &=
\langle P_{\hspace{-0.05cm}_{\mbox{\tiny $(23)$}}}(T), H_1
\circ_{\hspace{-0.05cm}_{\mbox{\tiny $(12)(3)$}}} H_2
\circ_{\hspace{-0.05cm}_{\mbox{\tiny $(12)(3)$}}} H_3 \rangle,
\end{align*}
and
\begin{align*}
\langle P_{\hspace{-0.05cm}_{\mbox{\tiny $(13)$}}}(T), H_1
\circ_{\hspace{-0.05cm}_{\mbox{\tiny $(13)(2)$}}} H_2
\circ_{\hspace{-0.05cm}_{\mbox{\tiny $(13)(2)$}}} H_3 \rangle &=
\langle P_{\hspace{-0.05cm}_{\mbox{\tiny $(13)$}}}(T), H_1
\circ_{\hspace{-0.05cm}_{\mbox{\tiny $(1)(2)(3)$}}} H_2
\circ_{\hspace{-0.05cm}_{\mbox{\tiny $(1)(2)(3)$}}} H_3 \rangle, \\
\langle P_{\hspace{-0.05cm}_{\mbox{\tiny $(23)$}}}(T), H_1
\circ_{\hspace{-0.05cm}_{\mbox{\tiny $(1)(23)$}}} H_2
\circ_{\hspace{-0.05cm}_{\mbox{\tiny $(1)(23)$}}} H_3 \rangle &=
\langle P_{\hspace{-0.05cm}_{\mbox{\tiny $(23)$}}}(T), H_1
\circ_{\hspace{-0.05cm}_{\mbox{\tiny $(1)(2)(3)$}}} H_2
\circ_{\hspace{-0.05cm}_{\mbox{\tiny $(1)(2)(3)$}}} H_3 \rangle.
\end{align*}
Finally, for any two permutations $\sigma_1, \sigma_2$ on $\N_3$
we have
\begin{align*}
\langle P_{\hspace{-0.05cm}_{\mbox{\tiny $(123)$}}}(T), H_1
\circ_{\sigma_1} H_2 \circ_{\sigma_1} H_3 \rangle &= \langle
P_{\hspace{-0.05cm}_{\mbox{\tiny $(123)$}}}(T), H_1
\circ_{\sigma_2} H_2 \circ_{\sigma_2} H_3 \rangle.
\end{align*}
\end{corollary}

\begin{example} \rm
\label{ex:1}
 In this example we demonstrate that
Formula~(\ref{deriv:1}) for the first derivative of a spectral
function, at $X$, can be rewritten in a different form. Let $X=V
(\Diag \lambda(X)) V^T$ and $\tilde{E} = V^T EV$, where $E$ is a
symmetric matrix. Using the definitions and notation in the
previous subsection we have:
\begin{align*}
\nabla (f \circ \lambda)(X)[E] &= \langle V \big( \Diag
\nabla f(\mu) \big) V^T, E \rangle \\
&= \langle \nabla f(\mu), \diag \tilde{E} \rangle \\
&= \langle \nabla f(\mu), \circ_{\hspace{-0.05cm}_{\mbox{\tiny
$(1)$}}} \tilde{E} \rangle.
\end{align*}
\end{example}

\begin{example} \rm
\label{ex:2} \label{Hes:exmpl} Let $X$ be a symmetric matrix with
ordered spectral decomposition $X=V (\Diag \lambda(X)) V^T$. Take
two symmetric matrices $E_{1}$ and $E_{2}$ and let $\tilde{E_i} =
V^T E_iV$ for $i=1,2$. As we saw in the examples in
Section~\ref{GenHadProd} we have:
$$
E_{1} \circ_{\hspace{-0.05cm}_{\mbox{\tiny $(1)(2)$}}} E_{2} =
(\diag E_{1})  (\diag E_{2})^T \,\, \mbox{ and } \,\, E_{1}
\circ_{\hspace{-0.05cm}_{\mbox{\tiny $(12)$}}} E_{2} = E_{1} \circ
E_{2}.,
$$
Then Formula~(\ref{deriv:2}) for the Hessian of the spectral
function $f \circ \lambda$ becomes:
\begin{align*}
\nabla^2 (f \circ \lambda)(X)[E_1,E_2] &= \nabla^2
f(\lambda(X))[\diag \tilde{E}_{1},\diag \tilde{E}_{2}] + \langle
\mathcal{A}(\lambda(X)), \tilde{E}_{1} \circ \tilde{E}_{2} \rangle\\
&=\langle \nabla^2 f(\lambda(X)), \tilde{E}_{1}
\circ_{\hspace{-0.05cm}_{\mbox{\tiny $(1)(2)$}}} \tilde{E}_{2}
\rangle + \langle \mathcal{A}(\lambda(X)), \tilde{E}_{1}
\circ_{\hspace{-0.05cm}_{\mbox{\tiny $(12)$}}} \tilde{E}_{2}
\rangle.
\end{align*}
\end{example}

All these examples support the following conjecture, which
describes the structure of the higher-order derivatives of
spectral functions.

\begin{conjecture} \rm
\label{conj:1} The spectral function $f \circ \lambda$ is $k$
times (continuously) differentiable at $X$ if and only of $f(x)$
is $k$ times (continuously) differentiable at the vector
$\lambda(X)$. Moreover, there are $k$-tensor valued maps
$\mathcal{A}_{\sigma} : \R^n \rightarrow T^{k,n}$, $\sigma \in
P^k$, such that for any symmetric matrices $E_1$,...,$E_k$ we have
$$
\nabla^k (f \circ \lambda)(X)[E_1,...,E_k] =  \sum_{\sigma \in
P^k} \langle \mathcal{A}_{\sigma}(\lambda(X)), \tilde{E}_1
\circ_{\sigma} \cdots \circ_{\sigma} \tilde{E}_k \rangle,
$$
where $X=V (\Diag \lambda(X)) V^T$ and $\tilde{E}_i = V^T E_iV$,
for $i=1,..,k$.
\end{conjecture}

In \cite{Sendov2003c}  we will show that this conjecture holds for
the derivatives of any spectral function at a symmetric matrix $X$
with distinct eigenvalues, as well as for the derivatives of
separable spectral functions at an arbitrary symmetric matrix.
(Separable spectral functions are those arising from symmetric
functions $f(x) = g(x_1) + \cdots + g(x_n)$ for some function $g$
on a scalar argument.) There we also describe how to compute the
operators $\mathcal{A}_{\sigma}$ for every $\sigma$ in $P^k$.

There is one major draw-back of the conjectured formula above. On
the left hand-side we have the the $k$-th derivative of the
spectral function evaluated at the matrices $E_1,...,E_k$ while on
the right-hand side these matrices are ``jumbled'' with the
orthogonal matrix $V$ into the $\sigma$-Hadamard products
$\tilde{E}_1 \circ_{\sigma} \cdots \circ_{\sigma} \tilde{E}_k$.
This is the problem that we address in the next section.

\section{The $\Diag^{\sigma}$ operator}
\label{STI}

Recall that the adjoint of the linear operator $\Diag : \R^n
\rightarrow M^n$ is the operator $\diag : M^n \rightarrow \R^n$.
That is, we have the identity
\begin{equation}
\label{sam} \langle \Diag x, H \rangle = \langle x, \diag H
\rangle,
\end{equation}
for any vector $x$ and any matrix $H$. It is also easy to verify
that for any vector $x$, matrix $H$, and orthogonal matrix $U$ we
have
\begin{equation}
\label{gen} \langle U (\Diag x) U^T , H \rangle = \langle x, \diag
(U^THU) \rangle.
\end{equation}
Vector $x$ can be viewed as a $1$-tensor on $\R^n$ given through
the linear isometry $x \rightarrow \langle x, \cdot \rangle$ and
similarly $\Diag x$ can be viewed as a $2$-tensor. In this section
we will generalize Equations~(\ref{sam}) and (\ref{gen}) for an
arbitrary $k$-tensor in place of $x$ and arbitrary
$\sigma$-Hadamard product in place of $\diag$.

Let $T$ be an arbitrary $k$-tensor on $\R^n$ and let $\sigma$ be a
permutation on $\N_k$. We define $\Diag^{\sigma} T$ to be a
$2k$-tensor on $\R^n$ in the following way
\begin{align*}
(\Diag^{\sigma} T)^{\substack{i_1...i_k \\ j_1...j_k}} = \left\{
\begin{array}{ll}
T^{i_1...i_k}, & \mbox{ if } i_s = j_{\sigma(s)}, \forall s = 1,...,k, \\
0, & \mbox{ otherwise. }
\end{array}
\right.
\end{align*}

When $k=1$ and $\sigma$ is the only choice from $P^1$, namely
$\sigma = (1)$, then this definition coincides with the definition
of the $\Diag$ operator in Equation~(\ref{sam}). Equivalent way to
define $\Diag^{\sigma} T$ that is useful for calculations is:
$$
(\Diag^{\sigma} T)^{\substack{i_1...i_k \\ j_1...j_k}} =
T^{i_1...i_k} \delta_{i_1j_{\sigma(1)}} \cdots
\delta_{i_kj_{\sigma(k)}}.
$$

We now define an action of the group, $O^n$, of all $n \times n$
orthogonal matrices on the space of all $k$-tensors on $\R^n$. For
any $k$-tensor $T$, and $U \in O^n$ this action will be denoted by
$UTU^T$, and defined by:
\begin{equation}
\label{kur0} (UTU^T)^{i_1...i_k} = \sum_{p_1 = 1 }^{n} \cdots
\sum_{p_k =1}^{n} \Big( T^{p_1...p_k} U^{i_{1}p_{1}} \cdots
U^{i_{k}p_{k}} \Big).
\end{equation}
In the case $k=1$, when $T$ is viewed as an $n$-dimensional
vector, this is exactly the action of the orthogonal group on
$\R^n$:
$$
(UTU^T)^{i_1} \equiv (UT)^{i_1} = \sum_{p_1=1}^n
U^{i_1p_1}T^{p_1}.
$$
In the case $k=2$ the definition coincides with the conjugate
action of the orthogonal group on the set of all $n \times n$
square matrices:
$$
(UTU^T)^{ij} = \sum_{p,q = 1}^{n,n} T^{pq}U^{ip}U^{jq},
$$
hence the use of the same notation for the general action $UTU^T$.
For future reference we state the formula of the action in the
case when the tensor is of even order. That is, if $T$ is a
$2k$-tensor, then
\begin{equation}
\label{kur} (UTU^T)^{\substack{i_1...i_k \\ j_1...j_k}} =
\sum_{p_1,q_1 = 1 }^{n,n} \cdots \sum_{p_k, q_k =1}^{n,n} \Big(
T^{\substack{p_1...p_k
\\ q_1...q_k}} \prod_{\nu =1 }^k
U^{i_{\nu}p_{\nu}}U^{j_{\nu}q_{\nu}} \Big).
\end{equation}

Let $P$ be an $n \times n$ permutation matrix and $\sigma$ its
corresponding permutation on $\N_n$, that is, $P^Te^i =
e^{\sigma(i)}$ for all $i=1,...,n$, where $\{e^i \, | \,
i=1,...,n\}$ is the standard basis in $\R^n$. The action of $P$ on
the tensors will be given by:
\begin{align*}
(PTP^T)^{i_1...i_k} &= \sum_{p_1 = 1 }^{n} \cdots \sum_{p_k
=1}^{n} \Big( T^{p_1...p_k} \prod_{\nu =1 }^k P^{i_{\nu}p_{\nu}}
\Big) = T^{\sigma(i_1)...\sigma(i_k)}.
\end{align*}
That is, the conjugate action of a permutation matrix on a
$k$-tensor is what one expects it to be. We have the following
immediate observation.
\begin{lemma}
For any permutation $\mu$ on $\N_k$, any permutation matrix $P$ in
$P^n$ and any $k$-tensor $T$ on $\R^n$, we have
$$
P(\Diag^{\mu} T)P^T = \Diag^{\mu} (PTP^T).
$$
\end{lemma}

\begin{proof}
Let $\sigma$ be the permutation on $\N_n$ corresponding to $P$.
Fix any multi index $(\substack{i_1...i_k \\ j_1...j_k})$. We
begin calculating the right-hand side entry corresponding to that
index. In the third equality below, we use the fact that $\sigma$
is a one-to-one map.
\begin{align*}
\big(P(\Diag^{\mu} T)P^T \big)^{\substack{i_1...i_k \\ j_1...j_k}}
&= \big(\Diag^{\mu} T \big)^{\substack{\sigma(i_1)...\sigma(i_k) \\
\sigma(j_1)...\sigma(j_k)}} \\
&= T^{\sigma(i_1)...\sigma(i_k)}
\delta_{\sigma(i_1)\sigma(j_{\mu(1)})} \cdots
\delta_{\sigma(i_k)\sigma(j_{\mu(k)})} \\
&= T^{\sigma(i_1)...\sigma(i_k)} \delta_{i_1j_{\mu(1)}} \cdots
\delta_{i_kj_{\mu(k)}} \\
&= (PTP^T)^{i_1...i_k} \delta_{i_1j_{\mu(1)}} \cdots
\delta_{i_kj_{\mu(k)}} \\
&= \big(\Diag^{\mu} (PTP^T) \big)^{\substack{i_1...i_k \\
j_1...j_k}}. \\[-1.5cm]
\end{align*}
\hfill \qed
\end{proof}

A natural question to ask is whether the action defined above on
the space $T^{k,n}$ is associative.

\begin{lemma}
\label{propty:1} for any $k$-tensor, $T$, on $\R^n$ and any two
orthogonal matrices $U$, $V$ in $O^n$ we have
$$
V(UTU^T)V^T=(VU)T(VU)^T.
$$
\end{lemma}

\begin{proof}
The proof is a direct calculation using the definitions. On one
hand we have
\begin{align*}
(V&(UTU^T)V^T)^{i_1...i_k} = \sum_{p_1 =1}^{n} \cdots \sum_{p_k
=1}^{n} \Big( (UTU^T)^{p_1...p_k} \prod_{\nu =1 }^k
V^{i_{\nu}p_{\nu}} \Big)\\
&=\sum_{p_1 =1}^{n} \cdots \sum_{p_k =1}^{n} \Big( \Big(
\sum_{l_1=1 }^{n} \cdots \sum_{l_k =1}^{n} T^{l_1...l_k}
\prod_{\mu =1 }^k U^{p_{\mu}l_{\mu}} \Big) \prod_{\nu =1 }^k
V^{i_{\nu}p_{\nu}} \Big).
\end{align*}
On the other hand we have
\begin{align*}
((VU)T(VU)^T)^{i_1...i_k} = \sum_{l_1 = 1 }^{n} \cdots \sum_{l_k
=1}^{n} T^{l_1...l_k} \prod_{\mu =1 }^k (VU)^{i_{\mu}l_{\mu}}.
\end{align*}
Using that
$$
(VU)^{i_{\mu}l_{\mu}}=\sum_{p_{\mu} =1 }^{n} V^{i_{\mu} p_{\mu}
}U^{p_{\mu}  l_{\mu}},
$$
we get
\begin{align*}
\prod_{\mu =1 }^k (VU)^{i_{\mu}l_{\mu}} &= \prod_{\mu =1 }^k \Big(
\sum_{p_{\mu} =1 }^{n} V^{i_{\mu} p_{\mu} }U^{p_{\mu}
l_{\mu}}\Big) = \sum_{p_1 = 1 }^{n} \cdots \sum_{p_k =1}^{n} \Big(
\prod_{\mu =1 }^k V^{i_{\mu} p_{\mu} } U^{p_{\mu} l_{\mu}} \Big).
\end{align*}
Putting everything together and observing that we can exchange the
multiple sum $\sum_{p_1 = 1 }^{n} \cdots \sum_{p_k =1}^{n}$ with
the multiple sum  $\sum_{l_1 = 1 }^{n} \cdots \sum_{l_k =1}^{n}$
we finish the proof of the lemma. \hfill \qed
\end{proof}

Let us see now that conjugation with an orthogonal matrix is
orthogonal transformation on $T^{k,n}$. That is, it doesn't change
the norm of the tensor. In other words, if we define
$$
\|T\| := \sqrt{\langle T, T\rangle},
$$
then we have the following lemma.

\begin{lemma}
Let $T$ be a $k$-tensor on $\R^n$, and $U$ be any orthogonal
matrix in $O^n$, then
$$
\|UTU^T\| = \|T\|.
$$
\end{lemma}

\begin{proof}
Direct calculation of the quantity $\|UTU^T\|^2$ gives:
\begin{align*}
&\|UTU^T\|^2 = \langle UTU^T, UTU^T\rangle \\
&= \sum_{i_1,...,i_k=1}^{n,...,n} (UTU^T)^{i_1...i_k}(UTU^T)^{i_1...i_k} \\
&= \sum_{i_1,...,i_k=1}^{n,...,n}
\Big(\sum_{p_1,...,p_k=1}^{n,...,n}T^{p_1...p_k}U^{i_1p_1} \cdots
U^{i_kp_k} \Big) \Big(
\sum_{q_1,...,q_k=1}^{n,...,n}T^{q_1...q_k}U^{i_1q_1} \cdots
U^{i_kq_k}\Big) \\
&=\sum_{p_1,...,p_k=1}^{n,...,n}\sum_{q_1,...,q_k=1}^{n,...,n}
T^{p_1...p_k}T^{q_1...q_k}\Big(\sum_{i_1=1}^{n}U^{i_1p_1}U^{i_1q_1}\Big)
\cdots \Big(\sum_{i_k=1}^{n} U^{i_kp_k}U^{i_kq_k}\Big) \\
&= \sum_{p_1,...,p_k=1}^{n,...,n}\sum_{q_1,...,q_k=1}^{n,...,n}
T^{p_1...p_k}T^{q_1...q_k} \delta_{p_1q_1} \cdots \delta_{p_kq_k}\\
&= \sum_{p_1,...,p_k=1}^{n,...,n} (T^{p_1...p_k})^2 \\
&= \|T\|^2.  \\[-1.5cm]
\end{align*}
\hfill \qed
\end{proof}

After all these preparations, we can give the following
generalization to Equation~(\ref{gen}). (When, $k=1$ and $\sigma =
(1)$ we obtain exactly Equation~(\ref{gen}).)

\begin{theorem}
\label{jen3} For any $k$-tensor $T$, any matrices $H_1$,...,$H_k$,
any orthogonal matrix $U$, and any permutation $\sigma$ on $\N_k$
we have the identity
\begin{equation}
\label{gen:1} \langle T, \tilde{H}_1 \circ_{\sigma} \cdots
\circ_{\sigma} \tilde{H}_k \rangle = \big(U(\Diag^{\sigma}
T)U^T\big)[H_1,...,H_k],
\end{equation}
where $\tilde{H}_i = U^TH_iU$, for all $i=1,2,...,k$.
\end{theorem}

\begin{proof}
Since both sides are linear in each argument separately, it is
enough to show that the equality holds for $k$-tuples
$(H_{i_1j_1},...,H_{i_kj_k})$ of basic matrices.

Using Lemma~\ref{Hadam:arb} and the fact that $\tilde{H}_{ij}^{pq}
= U^{ip}U^{jq}$, we develop the left-hand side of
Equation~(\ref{gen:1}):
\begin{align*}
\langle T, \tilde{H}_{i_1j_1} \circ_{\sigma} \cdots \circ_{\sigma}
\tilde{H}_{i_kj_k} \rangle &= \sum_{p_1,...,p_k = 1}^{n,...,n}
 T^{p_1...p_k}\tilde{H}_{i_1j_1}^{p_1p_{\sigma^{-1}(1)}}\cdots
 \tilde{H}_{i_kj_k}^{p_kp_{\sigma^{-1}(k)}} \\
&=\sum_{p_1,...,p_k = 1}^{n,...,n}
 T^{p_1...p_k} U^{i_1p_1}U^{j_1p_{\sigma^{-1}(1)}} \cdots
 U^{i_kp_k}U^{j_kp_{\sigma^{-1}(k)}}.
\end{align*}
On the other hand, using the definitions we calculate that the
right-hand side is:
\begin{align*}
(U(\Diag^{\sigma} T)U^T&)[H_{i_1j_1},...,H_{i_kj_k}] = \\
&=\sum_{p_1, q_1 =1}^{n,n} \cdots \sum_{p_k, q_k =1}^{n,n}
\big( (U(\Diag^{\sigma} T)U^T)^{\substack{p_1...p_k \\
q_1...q_k}}H^{p_1q_1}_{i_1j_1} \cdots H^{p_kq_k}_{i_kj_k}\big) \\
&=(U(\Diag^{\sigma} T)U^T)^{\substack{i_1...i_k \\
j_1...j_k}} \\
&= \sum_{p_1,q_1 = 1 }^{n,n} \cdots \sum_{p_k, q_k =1}^{n,n} \Big(
(\Diag^{\sigma} T)^{\substack{p_1...p_k
\\ q_1...q_k}} \prod_{\nu =1 }^k
U^{i_{\nu}p_{\nu}}U^{j_{\nu}q_{\nu}} \Big) \\
&=\sum_{p_1 = 1 }^{n} \cdots \sum_{p_k =1}^{n} \Big( T^{p_1...p_k}
\prod_{\nu =1 }^k
U^{i_{\nu}p_{\nu}}U^{j_{\nu}p_{\sigma^{-1}(\nu)}} \Big).
\end{align*}
This shows that the both sides are equal. \hfill \qed
\end{proof}

If we take the orthogonal matrix $U$ to be the identity matrix we
obtain the following corollary.

\begin{corollary}
\label{corl3a} For any $k$-tensor $T$, any matrices
$H_1$,...,$H_k$, and any permutation $\sigma$ on $\N_k$, we have
the identity
\begin{equation}
\label{gen:1a} \langle T,  H_1 \circ_{\sigma}...\circ_{\sigma} H_k
\rangle = (\Diag^{\sigma} T)[H_1,...,H_k].
\end{equation}
\end{corollary}

If in Corollary~\ref{corl3a} we substitute the matrices
$H_1$,...,$H_k$ with $\tilde{H}_1$,...,$\tilde{H}_k$ and we use
Theorem~\ref{jen3}, we obtain the next result.

\begin{corollary}
\label{corl3b} For any $k$-tensor $T$, orthogonal matrix $U \in
O(n)$, permutation $\sigma$ on $\N_k$, and any matrices
$H_1$,...,$H_k$ we have the identity
\begin{equation}
\label{gen:1b} (\Diag^{\sigma} T)[\tilde{H}_1,...,\tilde{H}_k] =
\big(U(\Diag^{\sigma} T)U^T\big)[H_1,...,H_k].
\end{equation}
\end{corollary}

If in Corollary~\ref{corl3a} we take $\sigma$ to be the identity
permutation, then we get the next corollary, which generalizes
Equation~(\ref{sam}).

\begin{corollary}
\label{corl3c} For any $k$-tensor $T$, any matrices
$H_1$,...,$H_k$ we have the identity
\begin{equation}
\label{gen:1c} T[\diag H_1,...,\diag H_k] =
(\Diag^{(\mbox{\scriptsize \rm id})} T)[H_1,...,H_k].
\end{equation}
\end{corollary}

We conclude this section with a second look at the first two
derivatives of spectral functions.

\begin{example} \rm
\label{ex:1b} As we saw in Example~\ref{ex:1}, the first
derivative of the spectral function $f \circ \lambda$ at the point
$X=V (\Diag \lambda(X)) V^T$, applied to the symmetric matrix $E$
is given by the formula
$$
\nabla (f \circ \lambda)(X)[E] = \langle \nabla f(\mu),
\circ_{\hspace{-0.05cm}_{\mbox{\tiny $(1)$}}} \tilde{E} \rangle,
$$
where $\tilde{E} = V^T EV$. This formula can be rewritten as
$$
\nabla (f \circ \lambda)(X)[E] = \langle V\big(\Diag \nabla f(\mu)
\big)V^T , E \rangle = V\big(\Diag \nabla f(\mu) \big)V^T[E].
$$
This was essentially the original form of this formula given in
Equation~(\ref{deriv:1}).
\end{example}

The usefulness of the new notation becomes more evident below.

\begin{example} \rm
\label{ex:2b} Let $X$ be a symmetric matrix with ordered spectral
decomposition $X=V (\Diag \lambda(X)) V^T$. Take two symmetric
matrices $E_{1}$ and $E_{2}$ and let $\tilde{E_i} = V^T E_iV$ for
$i=1,2$. As we saw in Example~\ref{ex:1}, the Hessian of the
spectral function $f \circ \lambda$ at the point $X=V (\Diag
\lambda(X)) V^T$, applied to the symmetric matrices $E_1$ and
$E_2$ is given by the formula
\begin{align*}
\nabla^2 (f \circ \lambda)(X)[E_1,E_2] &=\langle \nabla^2
f(\lambda(X)), \tilde{E}_{1} \circ_{\hspace{-0.05cm}_{\mbox{\tiny
$(1)(2)$}}} \tilde{E}_{2} \rangle + \langle
\mathcal{A}(\lambda(X)), \tilde{E}_{1}
\circ_{\hspace{-0.05cm}_{\mbox{\tiny $(12)$}}} \tilde{E}_{2}
\rangle.
\end{align*}
With the notation introduced in this section we can rewrite it as
\begin{align*}
\nabla^2 (f \circ \lambda)(X)[E_1,E_2] &= \big(V\big(
\Diag^{(1)(2)} \nabla^2 f(\lambda(X)) \big)V^T \big) [ E_{1},
E_{2} ] \\
& \hspace{3.2cm} + \big( V \big( \Diag^{(12)}
\mathcal{A}(\lambda(X)) \big) V^T \big)[ E_{1}, E_{2}].
\end{align*}
Or, in other words
$$
\nabla^2 (f \circ \lambda)(X) = V\big( \Diag^{(1)(2)} \nabla^2
f(\lambda(X)) + \Diag^{(12)} \mathcal{A}(\lambda(X)) \big) V^T.
$$
\end{example}

Finally, we express Conjecture~\ref{conj:1} in the new language.

\begin{conjecture} \rm
\label{conj2} The spectral function $f \circ \lambda$ is $k$ times
(continuously) differentiable at $X$ if, and only if, $f(x)$ is
$k$ times (continuously) differentiable at the vector
$\lambda(X)$. Moreover, there are $k$-tensor valued maps
$\mathcal{A}_{\sigma} : \R^n \rightarrow T^{k,n}$, $\sigma \in
P^k$, such that
\begin{equation}
\label{conj:eqn} \nabla^k (f \circ \lambda)(X) = V \Big(
\sum_{\sigma \in P^k} \Diag^{\sigma}
\mathcal{A}_{\sigma}(\lambda(X)) \Big) V^T,
\end{equation}
where $X=V (\Diag \lambda(X)) V^T$.
\end{conjecture}

In \cite{Sendov2003c}  we will show that this conjecture holds for
the derivatives of any spectral function at a symmetric matrix $X$
with distinct eigenvalues, as well as for the derivatives of
separable spectral functions at an arbitrary symmetric matrix.
(Separable spectral functions are those arising from symmetric
functions $f(x) = g(x_1) + \cdots + g(x_n)$ for some function $g$
on a scalar argument.) There we also describe how,  for every
$\sigma$ in $P^k$, to compute the operators
$\mathcal{A}_{\sigma}$, that depend only on the symmetric function
$f$.

\section{Comments on Conjecture~\ref{conj2}}
\label{last-sect}

In this section we show that once Conjecture~\ref{conj2} is
established for $k=1$, then for $k \ge 2$ it is enough to prove it
only in the case when the $X = \Diag x$ for some $x \in \R^n$ with
$x_1 \ge \cdots \ge x_n$. We begin with a simple lemma. For
brevity, given a $k$-tensor, $T$, on $M^n$ by $T[H]$ we denote the
$(k-1)$-tensor $T[\cdot,...,H]$.

\begin{lemma}
\label{simplelem} Let $T$ be any $2k$-tensor on $R^n$, $U \in
O^n$, and let $H$ be any matrix. Then, the following identity
holds.
$$
U(T[\tilde{H}])U^T = (UTU^T)[H],
$$
where $\tilde{H}=U^THU$.
\end{lemma}

\begin{proof}
Since both sides are linear with respect to $H$, it is enough to
prove the identity only for basic matrices $H_{i_kj_k}$. By the
definition of conjugation, and using the fact that
$\tilde{H}_{i_kj_k}^{pq} = U^{i_kp}U^{j_{k}q}$ we obtain
\begin{align*}
\big(U(T[\tilde{H}_{i_kj_k}])&U^T\big)^{\substack{i_1...i_{k-1} \\
j_1...j_{k-1}}} \\
&=\sum_{\substack{p_s, q_s = 1 \\
s=1,...,k-1}}^{n,...,n} (T[\tilde{H}_{i_kj_k}])^{\substack{p_1...p_{k-1} \\
q_1...q_{k-1}}} U^{i_1p_1}U^{j_1q_1} \cdots
U^{i_{k-1}p_{k-1}}U^{j_{k-1}q_{k-1}} \\
&= \sum_{\substack{p_s, q_s = 1 \\
s=1,...,k}}^{n,...,n} T^{\substack{p_1...p_{k} \\
q_1...q_{k}}} U^{i_1p_1}U^{j_1q_1} \cdots
U^{i_{k}p_{k}}U^{j_{k}q_{k}} \\
&= (UTU^T)^{\substack{i_1...i_{k} \\
j_1...j_{k}}} \\
&= \big((UTU^T)[H_{i_kj_k}] \big)^{\substack{i_1...i_{k-1} \\
j_1...j_{k-1}}}. \\[-1.4cm]
\end{align*}
\hfill \qed
\end{proof}

Suppose that Conjecture~\ref{conj2} holds for all derivatives of
order less than $k$ and for the $k$-th derivative it holds only
for ordered diagonal matrices. We will show that the conjecture
holds for the $k$-th derivative at an arbitrary matrix. Indeed,
let $X=V (\Diag \lambda(X)) V^T$, let $E$ be arbitrary symmetric
matrix and denote $\tilde{E}=V^TEV$. Then
\begin{align*}
&\nabla^{k-1} F(X+E) = \nabla^{k-1} F\big(V
(\Diag \lambda(X) + \tilde{E}) V^T \big) \\
&= V\big( \nabla^{k-1} F(\Diag \lambda(X)+\tilde{E}) \big)V^T  \\
&= V\big( \nabla^{k-1} F(\Diag \lambda(X)) \big)V^T + V\big(
\nabla^{k} F(\Diag \lambda(X))[\tilde{E}] \big)V^T + o(\|E\|) \\
&= \nabla^{k-1} F(X) + \big(V(\nabla^{k} F(\Diag \lambda(X)))V^T
\big)[E] + o(\|E\|).
\end{align*}
This shows that $\nabla^{k-1}F$ is differentiable at $X$ and that
$V(\nabla^{k} F(\Diag \lambda(X)))V^T$ is the $k$-th derivative of
$F$ at $X$.

\begin{proposition}
Suppose the $k$-th derivative of the spectral function $F=f \circ
\lambda$ is given by Equation~(\ref{conj:eqn}) for all $X$. If for
every $\sigma \in P^k$ the tensor valued map $x \in \R^n
\rightarrow \mathcal{A}_{\sigma}(x) \in T^{k,n}$ is continuous,
then $\nabla^kF(X)$ is continuous in $X$, in other words $F \in
C^k$.
\end{proposition}

\begin{proof}
Suppose that there is a sequence of symmetric matrices $X_m$
approaching $X$ and an $\epsilon > 0$ such that
$$
\|\nabla^kF(X_m) - \nabla^kF(X) \| > \epsilon, \mbox{ for all } m.
$$
Let $X_m=V_m (\Diag \lambda(X_m)) V_m^T$ and suppose without loss
of generality that the orthogonal $V_m$ approaches $V$.
(Otherwise, take a subsequence.) Clearly, we have $X=V (\Diag
\lambda(X)) V^T$, and by continuity of eigenvalues $\lambda(X_m)$
approaches $\lambda(X)$. Using the formula for the $k$-th
derivative and the continuity of the tensorial maps, the
contradiction follows. \hfill \qed
\end{proof}

\section{Equivalence relations on $\N_n$}

Suppose that $\sim$ is an equivalence relation on the integers
$\N_n$ and denote by $I_1$, $I_2$,...,$I_r$ the equivalence
classes determined by $\sim$. The equivalence classes will be also
called {\it blocks}. One may assume that the blocks are numbered
so that $I_1$ contains the integer $1$, $I_2$ contains that
smallest integer not in $I_1$, $I_3$ contains the smallest integer
not in $I_1 \cup I_2$, and so on.

In this short section, we will be interested in tensors having the
following structure.

\begin{definition} \rm
We say that a tensor $T \in T^{k,n}$ is {\it block-constant} (with
respect to the equivalence relation $\sim$) if
$$
T^{i_1...i_k} = T^{j_1...j_k}, \, \mbox{ whenever } \, i_s \sim
j_s \, \mbox{ for all } s=1,2,...,k.
$$
\end{definition}

Let $\mu$ be an arbitrary but fixed permutation in $P^k$.  We
introduce the linear operator $\tilde{P}_{\mu}$ on the space
$T^{k,n}$, generalizing the operator $P_{\mu}$ defined in
Section~\ref{sect:part-ord}. The definition is element-wise, as
follows:
$$
\big(\tilde{P}_{\mu}(T) \big)^{i_1...i_k} := \left\{
\begin{array}{ll}
T^{i_1...i_k}, & \mbox{ if } i_s \sim i_{\mu(s)} \, \forall s \in
\N_k, \\
0, & \mbox{ otherwise. }
\end{array}
\right.
$$

Clearly, when the equivalence relation $\sim$ is such that $i \sim
j$ if, and only if, $i=j$, then $\tilde{P}_{\mu}$  becomes equal
to the previously defined $P_{\mu}$. We would like to conclude
this work with a generalization of Theorem~\ref{trans:thm}.

\begin{theorem}
\label{trans-block:thm} Let $\sigma_1$, $\sigma_2$, and $\mu$ be
three permutations in $P^k$. Then for any block-constant matrices
$H_1$,...,$H_k$, and any tensor $T$ in $T^{k,n}$ we have the
identity:
\begin{equation*}
\langle \tilde{P}_{\mu}(T), H_1 \circ_{\sigma_1} \cdots
\circ_{\sigma_1} H_k \rangle = \langle \tilde{P}_{\mu}(T), H_1
\circ_{\sigma_2} \cdots \circ_{\sigma_2} H_k \rangle
\end{equation*}
if, and only if, $\mu \preceq \sigma_2^{-1} \circ \sigma_1$.
\end{theorem}

\begin{proof}
The proof is completely analogous to the one of
Theorem~\ref{trans:thm}. Consider a basis for the space of block
constant matrices  $\{ H_{pq} \, : \, 1 \le p,q \le n \}$ such
that $H_{pq}^{ij}$ is equal to one, if $i \sim j$, and zero
otherwise. Then all we have to change in the proof of
Theorem~\ref{trans:thm} is the ``$=$'' signs between indexes with
``$\sim$'' signs and all ``$\not=$'' signs with ``$\nsim$''.
\hfill \qed
\end{proof}


%

\end{document}